\documentclass[final,twoside,a4paper,12pt]{amsart}

\usepackage[notref]{showkeys}

\newcounter{exm}
\newcounter{manev}

\newcommand{\HH}{\mathfrak{H}}
\newcommand{\nep}{{1-\varepsilon}}
\newcommand{\dep}{{2\varepsilon}}
\newcommand{\KK}{\mathfrak{K}}
\newcommand{\kk}[1]{\mathfrak{H}_{#1}}

\newcommand{\mf}[1]{\mathfrak{#1}}
\newcommand{\eps}{\varepsilon}
\newcommand{\BH}{\mathcal{B}(\mf{H})}

\newtheorem{prop}{Proposition}[section]
\newtheorem{cor}[prop]{Corollaire}
\newtheorem{deff}{D\'efinition}[section]
\newtheorem{deffs}{D\'efinition}
\newtheorem{lemme}[prop]{Lemme}
\newtheorem{thm}[prop]{Th\'eor\`eme}
\newtheorem{thms}{Th\'eor\`eme}
\newtheorem{rem}[prop]{Remarque}

\newenvironment{dem}{\textbf{D\'emonstration. }}{\qed\vspace{0.5cm}}

\title[Bimodules repr\'esentables]{Bimodules norm\'es repr\'esentables sur des espaces hilbertiens}
\author{Ciprian Pop}

\begin{document}
\maketitle

\begin{section}{Introduction}
Soient $A$, $B$ deux $C^*$ alg\`ebres. Les articles~\cite{ces} et~\cite{er1}
contiennent une caract\'erisation abstraite des $A-B$ bimodules $V$
ayant une structure d'espace d'op\'erateurs telle qu'il existe un espace de
Hilbert $\mf{H}$, deux repr\'esentations fid\`eles 
$\pi:A\to\mathcal{B}(\mf{H})$, $\rho:B\to\mathcal{B}(\mf{H})$ et une 
application compl\`etement isom\'etrique $J:V\to\BH$ v\'erifiant
\begin{equation*}
  J(avb)=\pi(a)J(v)\rho(b),\quad\forall a\in A,b\in B,v\in V.
\end{equation*}
De tels bimodules interviennent dans la th\'eorie cohomologique des alg\`ebres
d'op\'erateurs. Ils sont souvent appel\'es $A-B$ bimodules 
d'o\-p\'e\-ra\-teurs.
Ici nous les appellerons plut\^ot $A-B$ bimodules $L^\infty$ matriciellement
norm\'es.

Etant donn\'e un $A-B$ bimodule norm\'e $V$, il est naturel de chercher des 
crit\`eres
assurant l'existence d'une structure de $A-B$ bimodule $L^\infty$ 
matriciellemnt norm\'e sur $V$ compatible avec la norme donn\'ee sur $V$,
c'est \`a dire telle que $\|\cdot\|_1=\|\cdot\|$. Il est clair que la condition
suivante, appel\'ee condition $(R)$, est n\'ec\'essaire:
\begin{equation*}
  \|a_1v_1b_1+a_2v_2b_2\|\le\|a_1a_1^*+a_2a_2^*\|^{1/2}
  \|b_1^*b_1+b_2^*b_2\|^{1/2}\max\{\|v_1\|,\|v_2\|\},
\end{equation*}
pour tous $a_1,a_2\in A$, $b_1,b_2\in B$ et $v_1,v_2\in V$.

Nous d\'emontrons la r\'eciproque. Plus pr\'ecis\'ement, nous d\'emontrons
qu'un $A-B$ bimodule norm\'e satisfait \`a la condition $(R)$ si et
seulement si il existe un espace de Hilbert $\mf{H}$, deux repr\'esentations
fid\`eles $\pi,\rho$ de $A$ et $B$ dans $\mf{H}$ et une application
isom\'etrique de $A-B$ bimodule $j:V\to\BH$ tels que
\begin{equation*}
  j(avb)=\pi(a)j(v)\rho(b),\quad\forall a\in A,b\in B,v\in V.
\end{equation*}
Nous dirons que de tels bimodules sont repr\'esentables ($(R)$ comme
$R$uan ou $R$epr\'esentable). Le cas $A=B=\mathbb{M}_n$ a \'et\'e trait\'e
par F. Lehner dans~\cite{leh1}. Le cas commutatif $A=B=C_0(X)$, o\`u les
actions \`a droite et \`a gauche co\^\i ncident est d\^u \`a K.H. Hofmann
\cite{KiRo}.

Notre travail est organis\'e comme suit.

On commence par d\'emontrer quelques r\'esultats techniques utiles 
pour la
suite. L'ingr\'edient principal est le lemme fondamental~\ref{pici}, qui 
sert \`a
caract\'eriser les formes lin\'eaires born\'ees sur les $A-B$ bimodules
repr\'esentables. 

Pour une repr\'esentation $(\pi,\mf{H})$ d'une $C^*$-alg\`ebre $C$ nous
introduisons la notion suivante de \emph{cyclicit\'e locale}.
\begin{deffs}
  Soit $C$ est une $C^*$-alg\`ebre. Une repr\'esentation $(\pi,H)$ de
  $C$ est appel\'ee \emph{localement cyclique} si il existe 
  $\mathcal{E}\subset H$ un sous-ensemble dense dans $H$ tel que pour tout
  $h_1,h_2,\ldots,h_n\in\mathcal{E}$ il existe $h\in\mathcal{E}$ tel que 
  $h_i\in\overline{\pi(C)h}$.
\end{deffs}

Cette propri\'et\'e a \'et\'e remarqu\'ee dans~\cite{sm1}, o\`u on consid\`ere
$A,B\subset\mathcal{B}(\mf{H})$ deux $C^*$ alg\`ebres et $V\subset\BH$ un
$A-B$ bimodule. Si $A\hookrightarrow\BH$ et $B\hookrightarrow\BH$ sont
localement cycliques, Smith d\'emontre que tout morphisme 
de $A-B$ bimodules $R:V\to\BH$ est compl\`etement born\'e, avec
\begin{equation*}
  \|R\|=\|R\|_{cb}.
\end{equation*}
On en d\'eduit (\`a l'aide du th\'eor\`eme de prolongement de Wittstock) que
$R$ se prolonge en un morphisme de $A-B$ bimodule $\tilde{R}$ de $\BH$
dans lui-m\^eme, avec $\|\tilde{R}\|=\|R\|$ (voir aussi~\cite{sue}).

Dans la section 2 nous d\'emontrons, sans utiliser le r\'esultat de Wittstock,
le th\'eor\`eme de prolongement suivant.
\begin{thms}
  Soient $(\pi,\mathfrak{K})$ et $(\rho,\mathfrak{H})$ deux repr\'esentations
  localement cycliques de $A$ et $B$ respectivement.
  Soient $W$ un $A-B$ bimodule v\'erifiant la propri\'et\'e $(R)$ 
  et $V\subset W$ un sous
  $A-B$ bimodule. Soit 
  \begin{equation*}
    R:V\to\mathcal{B}(\mathfrak{K},\mathfrak{H})
  \end{equation*}
  un morphisme contractif de $A-B$ bimodules.
  Alors il existe un morphisme contractif de $A-B$ bimodules
  \begin{equation*}
    \tilde{R}:W\to\mathcal{B}(\mathfrak{K},\mathfrak{H})
  \end{equation*}
  qui prolonge $R$.
\end{thms}
De plus, le th\'eor\`eme de Wittstock est une cons\'equence imm\'ediate
de ce r\'esultat, gr\^ace \`a une propri\'et\'e remarquable de l'alg\`ebre
$\mathcal{K}$ des op\'erateurs compacts en dimension infinie, \`a savoir que 
toutes les repr\'esentations de $\mathcal{K}$ sont localement cycliques.

Dans la section 3 nous \'etudions les repr\'esentations des $A-B$ bimodules
norm\'es, et nous caract\'erisons les $A-B$ bimodules repr\'esentables
comme \'etant exactement ceux qui poss\`edent la propri\'et\'e $(R)$. Si
$\kk{A''}$ et $\kk{B''}$ d\'esignent les espaces des repr\'esentations 
standards des alg\`ebres de von Neumann enveloppantes de $A''$ et $B''$ de
$A$ et $B$ respectivement, le $A-B$ bimodule $\mathcal{B}(\kk{B''},\kk{A''})$
joue un r\^ole important dans notre travail ainsi que 
$End_{A,B}(V,\mathcal{B}(\kk{B''},\kk{A''})$, espace norm\'e des morphismes
de $A-B$ bimodules de $V$ dans $\mathcal{B}(\kk{B''},\kk{A''})$. Celui-ci
g\'en\'eralise dans notre cadre le dual d'un espace norm\'e ordinaire. On 
l'appelle le dual standard du $A-B$ bimodule $V$.

Etant donn\'e un $A-B$ bimodule repr\'esentable, nous montrons dans la 
section 4 qu'il existe une plus petite structure $L^\infty$ matricielle
compatible, et nous donnons plusieurs caract\'erisations. Plus pr\'ecis\'ement,
on a
\begin{thms}
  Soit $(V,(\|\cdot\|_n)_n)$ un $A-B$ bimodule $L^\infty$
  matriciellement norm\'e.
  Les affirmations suivantes sont equivalentes:
  \begin{enumerate}
  \item[(a)] La structure matricielle de $V$ est minimale.
  \item[(b)] Il existe $\Omega$ un compact et un morphisme de $A-B$
    bimodules compl\`etement isom\'etrique $V\to
    C(\Omega,\mathcal{B}(\kk{B''}, \kk{A''}))$.
  \item[(c)] Tout morphisme isom\'etrique $j:V\to
    C(\Omega,\mathcal{B}(\kk{B''},\kk{A''}))$ est compl\`etement
    isom\'etrique.
  \item[(d)] Pour $[v_{kl}]\in\mathbb{M}_n(V)$ on a
    \begin{equation*}
    \|[v_{kl}]\|_n=\sup\{\|\sum_{k=1}^n\sum_{l=1}^na_kv_{kl}b_l\|\ :
    \ \|\sum_{k=1}^na_ka_k^*\|\le1,\|\sum_{l=1}^nb_l^*b_l\|\le1\}.      
    \end{equation*}
  \item[(e)] Pour tout $A-B$ bimodule $L^\infty$ matriciellement norm\'e, 
    tout morphisme contractif de $A-B$ bimodules $W\to V$ est
    compl\`etement contractif.
  \end{enumerate}
\end{thms}

Comme cons\'equence directe nous donnons une autre d\'emonstration du 
th\'eor\`eme de repr\'esentation de~\cite{er1}.

Finalement, dans la section 5 nous repr\'ennons les m\^emes questions, 
en ajoutant
des hypoth\`eses de dualit\'e. Nous d\'efinissons les $A-B$ bimodules 
repr\'esentables duaux (en enl\'evant la propri\'et\'e $(R)$
ce sont ceux introduits dans~\cite{KaRi}). Nous caract\'erisons ces objets,
avec ou sans hypoth\`eses de normalit\'e pour les actions de $A$ et $B$.

Ces r\'esultats seront utiles pour l'\'etude des produits tensoriels 
relatifs des bimodules repr\'esentables, actuellement en cours.
\end{section}

\begin{section}
{D\'efinitions et r\'esultats de base} Soit $V$ un
espace norm\'e \'equip\'e d'une structure de $A-B$ bimodule, o\`u $A$ et $B$
sont deux $C^*$-alg\`ebres. On suppose que $V$ est essentiel, c'est
\`a dire  que 
\begin{equation*} 
  \overline{AV}=\overline{VB}=V.
\end{equation*}

\begin{deff}
  \label{defrp}
  Si pout tout $a\in A$, $b\in B$, $v\in V$ on a
  \begin{equation*} 
    \|avb\|\le\|a\|\cdot\|v\|\cdot\|b\|,
  \end{equation*} 
  on dit que $V$ est un \emph{$A-B$ bimodule norm\'e}.
\end{deff}

\begin{deff}
  Si pour tout $a_1,a_2\in A$, $b_1,b_2\in B$ et $\|v_1\|,\|v_2\|\le1$
  on a
  $$
  \|a_1v_1b_1+a_2v_2b_2\|\le\|a_1a_1^*+a_2a_2^*\|^{1/2}\cdot
  \|b_1^*b_1+b_2^*b_2\|^{1/2}, \leqno (R)$$
  on dit que $V$ v\'erifie
  la propri\'et\'e $(R)$.
  
  Si pour tout $a_1,a_2\in A^+$, $b_1,b_2\in B^+$ et
  $\|v_1\|,\|v_2\|\le1$ on a
  $$
  \|a_1v_1b_1+a_2v_2b_2\|\le\|a_1^2+a_2^2\|^{1/2}\cdot\|b_1^2+b_2^2\|^{1/2},
  \leqno (R')$$
  on dit que $V$ v\'erifie la propri\'et\'e $(R')$.
\end{deff}

Evidemment la propri\'et\'e $(R)$ implique que $V$ est un $A-B$
bimodule norm\'e.

Il est facile de voir que si $V$ est un $A-B$ bimodule norm\'e alors
pour toute unit\'e approch\'ee $(a_\alpha)$ de $A$ et pour toute
unit\'e approch\'ee $(b_\beta)$ de $B$ on a
\begin{equation*}
  v=\lim_\alpha a_\alpha v=\lim_\beta vb_\beta,\quad\forall v\in V.
\end{equation*}

\begin{prop}
  \label{nonunital}
  Soit $V$ un $A-B$ bimodule norm\'e. Alors les propri\'et\'es $(R)$
  et $(R')$ sont \'equivalentes.
\end{prop}

\begin{dem}
  Evidemment $(R)$ implique $(R')$. Supposons que $V$ satisfait \`a la
  condition $(R)$.  Soient $a_1,a_2\in A$, $b_1,b_2\in B$ et
  $v_1,v_2\in V$ tels que $\|v_i\|\le 1$, $i=1,2$.  On va utiliser un
  r\'esultat de \cite{com1} ou ~\cite{ped}.  Si $C$ est une $C^*$-
  alg\`ebre et $c\in C$, pour tout $0<\varepsilon<1$ il existe un
  unique $c_\eps\in C$ tel que $c=c_\eps|c|^{1-\eps}$ et
  $|c_\eps|=|c|^\eps$.  Avec cette observation appliqu\'ee \`a
  $a_1^*$, $a_2^*$, $b_1$ et $b_2$, pour tout $\eps>0$ on a
  \begin{equation*}
    a_1v_1b_1+a_2v_2b_2=|a_1^*|^{(1-\eps)}w_1|b_1|^{(1-\eps)}+
    |a_2^*|^{(1-\eps)}w_2|b_2|^{(1-\eps)},
  \end{equation*}
  o\`u $w_i=a'_iv_ib'_i$ avec $\|a'_i\|=\|a_i^*\|^\eps$ et
  $\|b'_i\|=\|b_i\|^\eps$, $i=1,2$. Soient
  \begin{eqnarray*}
    A_\eps&=&\||a_1^*|^{(2-2\eps)}+|a_2^*|^{(2-2\eps)}\|\\
    B_\eps&=&\||b_1|^{(2-2\eps)}+|b_2|^{(2-2\eps)}\|\\
    m_\eps&=&\max\{\|w_1\|,\|w_2\|\}.
  \end{eqnarray*}
  Par hypoth\`ese on a
  \begin{equation*}
    \|a_1v_1b_1+a_2v_2b_2\|\le A_\eps^{1/2}B_\eps^{1/2}m_\eps.
  \end{equation*}
  Maintenant il suffit de remarquer que
  \begin{eqnarray*}
    \lim_{\eps\to0}A_\eps&=&\|a_1a_1^*+a_2a_2^*\|\\
    \lim_{\eps\to0}B_\eps&=&\|b_1^*b_1+b_2^*b_2\|\\
    \overline{\lim_{\eps\to0}}\ m_\eps&\le&1,
  \end{eqnarray*}
  d'o\`u le r\'esultat de l'\'enonc\'e.
\end{dem}

Supposons que $V$ est un $A-B$ bimodule v\'erifiant la propri\'et\'e
$(R)$ et que $A$ n'a pas d'unit\'e. Alors il existe une action
naturelle sur $V$ de $\tilde{A}$ (la $C^*$-alg\`ebre d\'eduite de $A$
par adjonction d'une unit\'e), prolongeant l'action de $A$, d\'efinie
par
\begin{equation*}
  (a\oplus\lambda)v=av+\lambda v,\quad\forall a\oplus\lambda\in\tilde{A},v\in V.
\end{equation*}
On voit facilement que $V$ devient un $\tilde{A}-B$ bimodule qui
v\'erifie la propri\'et\'e $(R)$. La m\^eme chose est valable pour
l'action \`a droite de $B$. Donc, en g\'eneral, on pourra supposer que
$A$ et $B$ poss\`edent des unit\'es.

On aura besoin du r\'esultat technique suivant, qui appara\^\i t comme
cons\'equence directe de~\cite{com1} ou~\cite{ped}

  \begin{lemme}
    \label{tech}
    Soit $C$ une $C^*$-alg\`ebre. Soit $x_1,x_2,\ldots,x_n\in C$ et
    $0<\varepsilon<1$. Posons $x=(\sum x_k^*x_k)^{1/2}$. Alors il
    existe $c_1,c_2,\ldots,c_n\in C$ tels que
    $x_k=c_kx^{1-\varepsilon}$ pour $k=1,2,\ldots,n$ et
    \begin{equation*}
      \sum c_k^*c_k=x^{2\varepsilon}.
    \end{equation*}
  \end{lemme}

  Dans la suite nous supposons que tous les bimodules consid\'er\'es seront suppos\'es qu'ils
  satisfont la propri\'et\'e $(R)$.

  Un r\'esultat cl\'e est le suivant:
  \begin{lemme}[fondamental]
    \label{pici}
    Soit $F\in V^*$. Alors $\|F\|\le1$ si et seulement si il existent
    $\phi\in S(A)$, $\psi\in S(B)$ deux \'etats de $A$ et respectivement de $B$
    tels que
    \begin{equation*}
      |F(avb)|\le\phi(aa^*)^{1/2}\psi(b^*b)^{1/2}\|v\|,
    \end{equation*}
    pour tous $a\in A$, $b\in B$ et $v\in V$.
  \end{lemme}
  Des variantes se trouvent dans \cite[theorem C]{er2} et \cite{leh1}.
  
  \begin{dem}
    La condition est \'evidemment suffisante.
    
    Supposons que $\|F\|=1$.  Soit $C=A\oplus B$. Pour $c=a\oplus b\in
    C^+=A^+\oplus B^+$ posons
    \begin{equation*}
    G(c)=\sup_{\|v\|\le1}|F(a^{1/2}vb^{1/2})|.
  \end{equation*}
  Evidemment $0\le G(c)\le\|a^{1/2}\|\cdot
  \|b^{1/2}\|\le\max\{\|a\|,\|b\|\}=\|c\|$.  On va prouver que $G$ est
  concave sur $C^+$.  Soit $c_i=a_i\oplus b_i\in C^+$, $i=1,2$. Soit
  $v_1,v_2\in V_1$. Quitte a remplacer $v_1$ et $v_2$ par
  $\lambda_1v_1$ et $\lambda_2v_2$ avec $|\lambda_i|=1$, on peut
  supposer que $F(a_iv_ib_i)\ge0$, $i=1,2$.
  
  Soit $a=a_1+a_2$, $b=b_1+b_2$. Fixons $\varepsilon>0$.  Par le lemme
  \ref{tech} il existe $e_1,e_2\in A$ et $f_1,f_2\in B$ tels que
  \begin{equation*}
    \begin{matrix}
      a_1^{1/2}=a^{(1-\varepsilon)/2}e_1,
      &a_2^{1/2}=a^{(1-\varepsilon)/2}e_2,
      &e_1e_1^*+e_2e_2^*=a^\varepsilon,\\
      b_1^{1/2}=f_1b^{(1-\varepsilon)/2},
      &b_2^{1/2}=f_2b^{(1-\varepsilon)/2},
      &f_1^*f_1+f_2^*f_2=b^\varepsilon.
    \end{matrix}
  \end{equation*}
  Alors
  \begin{equation*}
    F(a_1^{1/2}v_1b_1^{1/2})+F(a_2^{1/2}v_2b_2^{1/2})=
    F(a^{(1-\varepsilon)/2}wb^{(1-\varepsilon)/2}),
  \end{equation*}
  o\`u
  \begin{equation*}
    w=e_1v_1f_1+e_2v_2f_2.
  \end{equation*}
  Mais alors $\|w\|\le\|e_1e_1^*+e_2e_2^*\|^{1/2}\cdot
  \|f_1^*f_1+f_2^*f_2\|^{1/2}=\|a\|^{\varepsilon/2}\|b\|^{\varepsilon/2}$.
  En passant \`a la limite pour $\varepsilon\to0$ on trouve
  \begin{equation*}
    F(a_1^{1/2}v_1b_1^{1/2})+F(a_2^{1/2}v_2b_2^{1/2})\le G(c_1+c_2),\quad
    \forall \|v_1\|,\|v_2\|\le1,
  \end{equation*}
  et donc
  \begin{equation*}
    G(c_1)+G(c_2)\le G(c_1+c_2).
  \end{equation*}
  Si $c_i\in C^+$, $i=1,2$, alors
  \begin{equation*}
    \|c_1\|-G(c_2)\ge\|c_2\|-\|c_1-c_2\|-G(c_2)\ge-\|c_1-c_2\|.
  \end{equation*}
  Soit $H:C_{sa}\to\mathbb{R}$ d\'efinie par
  \begin{equation*}
    H(h)=\inf\{\|c_1\|-G(c_2)\ :\ h=c_1-c_2\mbox{\ avec\ }c_1,c_2\in C^+\}.
  \end{equation*}
  Alors $H$ est une forme sous-lin\'eaire sur $C_{sa}$ et par le
  th\'eor\`eme de Hahn--Banach il existe une forme
  $\mathbb{R}$-lin\'eaire $\rho$ d\'efinie sur $C_{sa}$ telle que
  $\rho(h)\le H(h),\quad\forall h\in C_{sa}$.  En particulier, si
  $c\in C^+$ on aura $\rho(c)\le\|c\|$ et $-\rho(c)=\rho(-c)\le-G(c)$
  pour tout $c\in C^+$. Donc
  \begin{equation*}
    G(c)\le\rho(c)\le\|c\|,\quad\forall c\in C^+.
  \end{equation*}
  Alors l'extension complexe de $\rho$ \`a $C$ est un \'etat de $C$,
  car si $A$ et $B$ sont unif\`eres on a
  \begin{equation*}
    1=G(1)\le\rho(1)\le\|1\|=1,
  \end{equation*}
  sinon on a toujours
  \begin{equation*}
    \sup_{c\in C^+,\|c\|\le1}G(c)=1.
  \end{equation*}
  Maintenant on pose $s=\rho(1\oplus0)$, $t=\rho(0\oplus1)$. On
  remarque que $s+t=1$. En outre, $s,t>0$, parce que $F\not=0$. Soit
  $\phi(a)=s^{-1}\rho(a\oplus0)$, $\psi(b)=t^{-1}\rho(0\oplus b)$.
  Alors $\phi\in S(A)$, $\psi\in S(B)$ et
  \begin{equation*}
    |F(avb)|\le\rho(a^2\oplus b^2)=
    s\phi(a^2)+t\psi(b^2)\quad\forall a\in A^+,b\in B^+,\|v\|\le1,
  \end{equation*}
  et par le m\^eme raisonnement qu'\`a la Proposition~\ref{nonunital} on
  trouve que pour tout $a\in A$, $b\in B$ et $\|v\|\le1$ on a
  \begin{equation*}
    |F(avb)|\le
    s\phi(aa^*)+t\psi(b^*b).
  \end{equation*}
  Par un argument standard on trouve que
  \begin{equation*}
    |F(avb)|\le2\sqrt{st}
    \phi(aa^*)^{1/2}\psi(b^*b)^{1/2}\|v\|\quad\forall a\in A,b\in B,v\in V.
  \end{equation*}
  Mais $st\le1/4$, donc le lemme est demontr\'e. En fait, comme on a
  suppos\'e $\|F\|=1$ on a automatiquement $st=1/4$, donc $s=t=1/2$.
\end{dem}

\begin{lemme}
  \label{extension}
  Soit $W$ un $A-B$ bimodule v\'erifiant la condition $(R)$ et $V\subset W$ 
  stable par les 
  actions de $A$ et $B$. Soit $\phi\in S(A)$ et $\psi\in S(B)$ deux \'etats
  de $A$ et $B$ respectivement, et $F\in V^*$ une fonctionnelle lin\'eaire
  v\'erifiant, pour tout $a\in A$, $b\in B$ et $v\in V$,
  \begin{equation*}
    |F(avb)|\le\phi(aa^*)^{1/2}\psi(b^*b)^{1/2}\|v\|.
  \end{equation*}
  Alors il existe $\tilde{F}\in W^*$ une extension de $F$ telle que pour
  tout $a\in A$, $b\in B$ et $w\in W$ on a
  \begin{equation*}
    |F(awb)|\le\phi(aa^*)^{1/2}\psi(b^*b)^{1/2}\|w\|.
  \end{equation*}
\end{lemme}

\begin{dem}
  On peut supposer $A$ et $B$ unif\`eres.
  Consid\'erons $N:W\to\mathbb{R}_+$ l'application donn\'ee par
  \begin{equation*}
    N(w)=\inf_{w=aw'b}\frac{\phi(aa^*)+\psi(b^*b)}{2}\|w'\|.
  \end{equation*}
  En utilisant le lemme~\ref{tech} on obtient que
  \begin{equation*}
    N(w)=\inf\{\frac{\phi(a^2)+\psi(b^2)}{2}\|w'\|,
    \ w=aw'b,a\in A^+,b\in B^+\},
  \end{equation*}
  et finalement, pour tout $w\in W$ on a
  \begin{equation*}
    N(w)=\inf\{\frac{\phi(a^2)+\psi(b^2)}{2}\|w'\|,
    \ w=aw'b,a\in \mathcal{I}(A^+),b\in \mathcal{I}(B^+)\},
  \end{equation*}
  o\`u on $\mathcal{I}(A^+)$ et $\mathcal{I}(B^+)$ repr\'esentent les
  \'el\'ements positifs et inversibles de $A$ et $B$ respectivement.
  Avec les m\^emes arguments que dans la d\'emonstration de la 
  lemme~\ref{pici} on v\'erifie que $N$ est une semi-norme sur $W$, et
  en utilisant la derni\`ere forme de $N$, on obtient que
  \begin{equation*}
    |F(v)|\le N(v),\quad\forall v\in V.
  \end{equation*}
  Par le th\'eor\`eme de Hahn--Banach il existe
  une forme lin\'eaire $\tilde{F}\in W^*$ qui prolonge $F$ telle que
  \begin{equation*}
    |\tilde{F}(w)|\le N(w)\quad\forall w\in W.
  \end{equation*}
  La forme $\tilde{F}$ v\'erifie toutes les conditions requises.
\end{dem}

Si $(\pi,\mathfrak{H})$ et $(\rho,\mathfrak{K})$ sont deux
repr\'esentations de $A$ et $B$ respectivement,
$\mathcal{B}(\mathfrak{K},\mathfrak{H})$ sera toujours muni de sa structure de 
$A$--$B$ bimodule canonique, par les
actions naturelles
\begin{equation*}
  a.T.b=\pi(a)T\rho(b).
\end{equation*}
Evidemment $\mathcal{B}(\mathfrak{K},\mathfrak{H})$ poss\`ede le 
propri\'et\'e $(R)$.
  
  \begin{deff}
    Soit $C$ est une $C^*$-alg\`ebre. Une repr\'esentation $(\pi,H)$ de
    $C$ est appel\'ee \emph{localement cyclique} si il existe 
    $\mathcal{E}\subset H$ un sous-ensemble dense dans $H$ tel que pour tout
    $h_1,h_2,\ldots,h_n\in\mathcal{E}$ il existe $h\in\mathcal{E}$ avec
    $h_i\in\overline{\pi(C)h}$, $i=1,\ldots,n$. On dira que $\mathcal{E}$
    est un sousensemble localement cyclique de $H$.
  \end{deff}
  Les repr\'esentations localement cycliques jouent un r\^ole crucial dans ce
  travail.
  Dans~\cite{sm1} il est demontr\'e que, pour une alg\`ebre de 
  von Neumann $M$ agissant sur l'espace de Hilbert $H$, l'inclusion de
  $M$ dans $\mathcal{B}(H)$ est localement cyclique si 
  tout \'etat normal de $M'$ est un \'etat vectoriel.
  En particulier, si $M$ est une alg\`ebre de von Neumann, sa repr\'esentation
  standard est localement cyclique. 
  Pour une $C^*$ alg\`ebre $A$, sa repr\'esentation
  sur l'espace de Hilbert standard de son alg\`ebre de von Neumann
  enveloppante $A''$ est aussi localement cyclique. 

  \begin{thm}
    \label{ext}
    Soient $(\pi,\mathfrak{K})$ et $(\rho,\mathfrak{H})$ deux repr\'esentations
    localement cycliques de $A$ et $B$ respectivement.
    Soient $W$ un $A-B$ bimodule v\'erifiant la propri\'et\'e $(R)$ 
    et $V\subset W$ un sous
    $A-B$ bimodule. Soit 
    \begin{equation*}
      R:V\to\mathcal{B}(\mathfrak{K},\mathfrak{H})
    \end{equation*}
    un morphisme contractif de $A-B$ bimodules.
    Alors il existe un morphisme contractif de $A-B$ bimodules
    \begin{equation*}
      \tilde{R}:W\to\mathcal{B}(\mathfrak{K},\mathfrak{H})
    \end{equation*}
    qui prolonge $R$.
  \end{thm}
  \begin{dem}
    Soit $\mathcal{F}\subset\KK$ et $\mathcal{E}\subset\HH$ des sous ensembles
    denses dans $\KK$ et $\HH$ respectivement, satisfaisant aux conditions
    de localement cyclicit\'e.
    On d\'efinit sur $\mathcal{F}\times\mathcal{E}$ une relation d'ordre par
    \begin{equation*}
      (\eta,\xi)\prec(\eta',\xi')\Leftrightarrow 
      \eta\in\overline{\rho(B)\eta'},\ 
      \xi\in\overline{\pi(A)\xi'}.
    \end{equation*}
    Pour $(\eta,\xi)$ donn\'e, on consid\`ere la forme sur $V$
    \begin{equation*}
      F_{\eta,\xi}(v)=(R(v)\eta|\xi).
    \end{equation*}
    Par le lemme~\ref{extension} il existe une forme lin\'eaire 
    $\tilde{F}\in W^*$
    qui prolonge $F_{\eta,\xi}$ telle que
    \begin{equation*}
      |\tilde{F}(awb)|\le\|\pi(a)^*\xi\|\cdot\|\rho(b)\eta\|\cdot\|w\|,
      \quad\forall a\in A,b\in B,w\in W.
    \end{equation*}
    Donc pour chaque $w\in W$ on peut d\'efinir un  op\'erateur 
    \begin{equation*}
      R_{\eta,\xi}(w)\in
      \mathcal{B}(\overline{\rho(B)\eta},\overline{\pi(A)\xi})
    \end{equation*}
    par
    \begin{equation*}
      (R_{\eta,\xi}(w)\rho(b)\eta|\pi(a)\xi)=\tilde{F}_{\eta,\xi}(a^*wb).
    \end{equation*}
    
    Mais la relation d'ordre introduite auparavant est filtrante croissante, 
    car les deux repr\'esentations sont par hypoth\`ese localement cycliques.
    Soit $\mathcal{U}$ un ultrafiltre la contenant. Alors, 
    pour $w\in W$ fix\'e, on d\'efinit la forme sesquilin\'eaire
    \begin{equation*}
      \Phi_w(\eta_0,\xi_0)=\lim_\mathcal{U}(R_{\eta,\xi}(w)\eta_0|\xi_0)
    \end{equation*}
    sur $\mathfrak{K}\times\mathfrak{H}$. On montre assez
    facilement que $\Phi_w$ d\'efinit un op\'erateur
    \begin{equation*}
    \tilde{R}(w)\in\mathcal{B}(\KK,\HH)  
    \end{equation*}
    et que $w\mapsto\tilde{R}(w)$ v\'erifie les conditions
    requises.
  \end{dem}

Le th\'eor\`eme de prolongement ci-dessus n'est pas valable pour
un morphisme de $A-B$ bimodules arbitraire $R:V\to\mathcal{B}({\KK,\HH})$,
m\^eme sous les hypoth\`ese les plus simples $A=B=\mathbb{C}$.

Il semble que l'hypoth\`ese
de cyclicit\'e locale de $\pi$ et $\rho$ soit essentielle. Cependant, 
le th\'eor\`eme de Wittstock nous dit que, si $V$ est un espace d'op\'erateurs,
si $T:V\to\mathcal{B}(H)$ est une application \emph{compl\`etement born\'ee} et
 si $V\subset W$ compl\`etement isom\'etriquement, alors il existe un 
prolongement
$\tilde{T}:W\to\mathcal{B}(H)$ avec $\|\tilde{T}\|_{cb}=\|T\|_{cb}$. Ceci est
valable sans aucune hypoth\`ese sur l'espace de Hilbert $H$. La raison pour 
cela est que, si on note $\mathcal{K}$ l'alg\`ebre des op\'erateurs compacts
sur $l^2$, alors toutes les repr\'esentations nond\'eg\'en\'er\'ees de
$\mathcal{K}$ sont localement cycliques. En outre, si $C$ est une $C^*$ 
alg\`ebre quelconque et $(\pi,\mathfrak{H})$ est une repr\'esentation 
nond\'eg\'en\'er\'ee de $C$, alors 
$(\mathrm{id}\otimes\pi,l^2\otimes\mathfrak{H})$ est une repr\'esentation 
localement cyclique de $\mathcal{K}\otimes C$.

Dans la suite, pour une $W^*$-alg\`ebre $M$ on va noter par
$\mathfrak{H}_M$ l'espace de Hilbert de la repr\'esentation standard
de $M$ (cf.\cite{ha1}).  Le fait important est que tout \'etat normal
d'une alg\`ebre de von Neumann sous forme standard est un \'etat
vectoriel.

Si $C$ est une $C^*$ alg\`ebre, alors $C$ sera identifi\'ee avec l'image
dans $C''$ de l'inclusion canonique (o\`u $C''$ est son alg\`ebre de
von Neumann enveloppante).
  
  \begin{prop}
    \label{p2}
    Soit $F\in V^*_1$. Alors il existe $\eta\in\kk{B''}$ et
      $\xi\in\kk{A''}$ de norme $1$, $R:V\to\mathcal{B}(\kk{B''},\kk{A''})$ un
      morphisme contractif de $A$--$B$ bimodules tels que
      \begin{equation*}
        \label{main}
        F(\cdot)=(R(.)\eta|\xi).
      \end{equation*}
  \end{prop}

  \begin{dem}
    Par le lemme \ref{pici} il existe des
    \'etats $\varphi$ et $\psi$ de $A$ et $B$ respectivement, satisfaisant \`a
    \begin{equation*}
      \label{good}
      |F(avb)|\le\varphi(aa^*)^{1/2}\|v\|\psi(b^*b)^{1/2},\quad
      \forall v\in V,a\in A, b\in B.
    \end{equation*}
    Les extensions $\varphi^{**}$ et $\psi^{**}$ sur $A''$
    et $B''$ respectivement sont des \'etats normaux, qu'on note aussi
    $\varphi$ et $\psi$. Par \cite[Theorem 2.17]{ha1} ce sont des
    \'etats vectoriels, c'est \`a dire qu'il existe des vecteurs (de
    norme $1$) $\xi\in\kk{A''}$, $\eta\in\kk{B''}$ tels que
    $\varphi=\omega_\xi$, $\psi=\omega_\eta$.  On d\'efinit sur
    $B\eta\times A\xi$ la forme sesquilin\'eaire
    \begin{equation*}
      \label{forme}
      (b\eta,a\xi)\mapsto F(a^*vb),
    \end{equation*}
    qui \emph{est bien d\'efinie}, car
    \begin{equation*}
      \label{cucu}
      |F(a^*vb)|\le\|v\|\varphi(a^*a)^{1/2}\psi(b^*b)^{1/2}=
      \|v\|\cdot\|b\eta\|\cdot\|a\xi\|.
    \end{equation*}
    Pour chaque $v\in V$ il existe donc un unique op\'erateur 
    \begin{equation*}
      R(v)\in \mathcal{B}(\overline{B\eta},\overline{A\xi})
    \end{equation*}
    de norme inf\'erieure ou
    \'egale \`a $\|v\|$ tel que
    \begin{equation*}
      \label{cucu1}
      (R(v)b\eta|a\xi)=F(a^*vb),\quad\forall
      (a,b)\in A\times B.
    \end{equation*}
    Comme $\overline{B\eta}\subset\kk{B''}$ et 
    $\overline{A\xi}\subset\kk{A''}$ on peut 
    regarder $R(v)$ comme morphisme de $\kk{B''}$
    \`a $\kk{A''}$, not\'e toujours par $R(v)$.  On voit
    facilement que $R$ est un morphisme de $A$--$B$ bimodules, car
    \begin{eqnarray*}
      (R(avb)b_1\eta|a_1\xi)&=&
      F(a_1^*avbb_1)\\
      &=&(R(v)bb_1\eta|a^*a_1\xi)\\
      &=&(aR(v)bb_1\eta|
      a_1\xi),
    \end{eqnarray*}
    et donc
    \begin{equation*}
      R(avb)=aR(v)b\quad\forall a\in A,b\in B,v\in V.
    \end{equation*}
  \end{dem}

\end{section}

\begin{section}{Representations de $A-B$ bimodules}
  
  Soit $V$ un $A-B$ bimodule norm\'e.

  \begin{deff}
    Soit $(\pi,\mathfrak{H})$, $(\rho,\mathfrak{K})$ deux
    repr\'esentations fid\`eles de $A$ et $B$ respectivement et $\Omega$ 
    un espace compact. Soit
    $J:V\to C(\Omega, \mathcal{B}(\mathfrak{K},\mathfrak{H}))$ un
    morphisme de $A-B$ bimodules.  Un tel triplet $(J,\pi,\rho)$ est appel\'e
    une repr\'esentation du $A-B$ bimodule $V$.
  \end{deff}
  
  Pour simplifier les notations, lorsqu'il n'y a pas de risque de confusion
  on va identifier $A$ et $B$ avec leurs images dans
  $\mathcal{B}(\mathfrak{H})$ et $\mathcal{B}(\mathfrak{K})$ respectivement.
  
  \begin{deff}
    Soient $(\pi,\mathfrak{H})$, $(\rho,\mathfrak{K})$ des repr\'esentations de
    $A$ et $B$ respectivement.  On appelle \emph{dual associ\'e} \`a 
    $(\pi,\rho)$ l'espace
    \begin{equation*}
      End_{A,B}(V,\mathcal{B}(\KK,\HH))  ,
    \end{equation*}
    not\'e $V^\dag_{(\pi,\rho)}$. L'espace
    \begin{equation*}
      End_{A,B}(V,\mathcal{B}(\kk{B''},\kk{A''}))
    \end{equation*}
    est appel\'e 
    \emph{le dual standard} de $V$, not\'e $V^\dag$.
  \end{deff}

  Si $A=B=\mathbb{C}$ alors le dual standard co\"\i ncide avec le dual usuel.
  Si $A=B$ est une $C^*$ alg\`ebre, alors $V=A$ poss\`ede une structure
  de $A-A$ bimodule (par multiplications). En fait $A$ est toujours un
  $M(A)-M(A)$ bimodule v\'erifiant la propri\'et\'e $(R)$. 
  Si $(\pi,\HH)$ est une repr\'esentation de $A$, alors $V_{(\pi,\pi)}^\dag$
  s'identifie avec $\pi(A)'$.

  Si $V$ est un $A-B$ bimodule norm\'e, $V^\dag_{(\pi,\rho)}$ poss\`ede
  une structure naturelle de $\pi(A)'-\rho(B)'$ bimodule , et il v\'erifie
  la propri\'et\'e $(R)$ relativement \`a $\pi(A)'$ et $\rho(B)'$. Les
  actions de $\pi(A)'$ et $\rho(B)'$ sont d\'efinies par
  \begin{equation*}
    (a'.T.b')(v)=a'T(v)b',\quad\forall v\in V,a'\in\pi(A)',b'\in\pi(B)'.
  \end{equation*}

  Pour un $A-B$ bimodule norm\'e $V$ on va noter
  \begin{equation}
    \label{compact}
    \Omega_V=End_{A,B}(V,\mathcal{B}(\kk{B''},\kk{A''}))_1
  \end{equation}
  l'ensemble des morphismes de $A-B$ bimodules contractifs de $V$ \`a
  valeurs dans $\mathcal{B}(\kk{B''},\kk{A''})$. Autrement dit, $\Omega_V$
  est la boule unit\'e du dual standard de $V$.
  Il est facile de voir que $\Omega_V$ est compact pour la topologie
  de la convergence simple de $V$ dans $\mathcal{B}(\kk{B''},\kk{A''})$, 
  muni de la topologie $*$-faible.  

  On remarque que 
  \begin{equation*}
    J_V:V\to C(\Omega_V,\mathcal{B}(\kk{B''},\kk{A''}))
  \end{equation*}
  d\'efinie par $J_V(v)(R)=R(v)$ est une representation de $A-B$
  bimodules, appell\'ee \emph{repr\'esentation standard} de $V$.

  \begin{thm}
    \label{imp}
    Soit $V$ un $A-B$-bimodule poss\'edant la propri\'et\'e $(R)$.  
    Soit $v\in V$. Alors
    \begin{equation*}
      \|v\|=\sup\{\|R(v)\|\quad:\quad R\in \Omega_V\}.
    \end{equation*}
  \end{thm}
  \begin{dem}
    Soit $F\in V^*_1$ telle que
    $F(v)=1$. Alors par la proposition~\ref{p2} le r\'esultat est
    imm\'ediat.
  \end{dem}

  \begin{cor}
    \label{dual}
    Si $V$ est un $A-B$ bimodule poss\'edant la propri\'et\'e $(R)$, 
    la repr\'esentation
    standard $J_V:V\to C(\Omega_V,\mathcal{B}(\kk{B''},\kk{A''}))$ est
    une isom\'etrie
  \end{cor}
  
  On peut maintenant \'enoncer le r\'esultat suivant.
  \begin{thm}
    \label{repr_bimodules}
    Soit $V$ un $A-B$ bimodule norm\'e.  Alors les conditions
    suivantes sont \'equivalentes:
  \item{(1)} $V$ est un $A-B$ bimodule poss\'edant la propri\'et\'e $(R)$
    (ou $(R')$).
  \item{(2)} La repr\'esentation standard est une isom\'etrie.
  \item{(3)} Il existe une repr\'esentation isom\'etrique de $V$.
  \item{(4)} Il existe un espace de Hilbert $H$, deux
    repr\'esentations fid\`eles $\pi:A\to\mathcal{B}(H)$,
    $\rho:B\to\mathcal{B}(H)$
    et un morphisme isom\'etrique $J:V\to\mathcal{B}(H)$ tel que
    \begin{equation*}
      J(avb)=\pi(a)J(v)\rho(b)\quad,\forall a\in A,v\in V,b\in B.
    \end{equation*}
  \end{thm}

  \begin{dem}
    L'implication $(1)\Rightarrow(2)$ r\'esulte du th\'eor\`eme \ref{imp}.
    Les implications $(4)\Rightarrow(1)$ et $(2)\Rightarrow(3)$ sont
    imm\'ediates.
    
    Pour montrer que $(3)$ implique $(4)$, consid\'erons $J_0$ un
    morphisme isom\'etrique de $V$ dans $C(\Omega,\mathcal{B}(
    \mathfrak{K},\mathfrak{H}))$.  On prends simplement
    $H_0=\mathfrak{H}\oplus\mathfrak{K}$, les repr\'esentations
    $\pi\oplus0$ et $0\oplus\rho$ et $J(v)=\left[
      \begin{matrix}
        0&J_0(v)\\
        0&0
      \end{matrix}\right]$. 
    Donc $J$ v\'erifie les conditions requises et est \`a valeurs
    dans $C(\Omega,\mathcal{B}(H))$, qui est une $C^*$-alg\`ebre.
  \end{dem}

  \begin{deff}
    Un $A-B$ bimodule qui poss\`ede ces propri\'et\'es \'equivalentes sera
    dit $A-B$ bimodule repr\'esentable.
  \end{deff}

  On finit cette section par une caract\'erisation des inclusions de
  $A-B$ bimodules repr\'esentables.

  Par analogie avec le cas classique des espaces de Banach, pour
  tout morphisme de $A-B$ bimodules $T:V\to W$ on peut d\'efinir 
  le morphisme dual $T^\dag:W^\dag\to V^\dag$ par
  \begin{equation*}
    T^\dag(S)(v)=S(T(v)),\quad\forall S\in W^\dag,v\in V.
  \end{equation*}
  On obtient:

  \begin{cor}
    Soit $J:V\to W$ un morphisme de $A-B$ bimodules repr\'esentables.
    Alors les affirmations suivantes sont \`equivalentes:
    \begin{enumerate}
    \item[(i)] $J$ est une isom\'etrie.
    \item[(ii)] $J^\dag(\Omega_W)=\Omega_V$ (en particulier $J^\dag$ est
      surjective).
    \end{enumerate}
  \end{cor}
  \begin{dem}
    Pour d\'emontrer que $(i)$ implique $(ii)$,
    soit $R\in\Omega_W$. Alors
    \begin{equation*}
      J^\dag(R)(v)=R(J(v)),\quad\forall v\in V,
    \end{equation*}
    donc
    \begin{equation*}
      \|J^\dag(R)(v)\|\le\|J(v)\|=\|v\|,\quad\forall v\in V.
    \end{equation*}
    On a donc $J^\dag(\Omega_W)\subset\Omega_V$.
    
    Il reste \`a remarquer que, si on identifie $V$ \`a son image $J(V)$,
    pour tout $R\in\Omega_V$, par le theor\`eme~\ref{ext} il existe
    $\tilde{R}\in W^\dag$ le prolongeant, avec $\|\tilde{R}\|=\|R\|$.
    Autrement dit, 
    \begin{equation*}
      \tilde{R}(J(v))=R(v),\quad\forall v\in V,
    \end{equation*}
    d'o\`u $J^\dag(\tilde{R})=R$,
    donc $J^\dag(\Omega_W)=\Omega_V$.

    Supposons maintenant que $J^\dag(\Omega_W)=\Omega_V$. 
    Par le th\'eor\`eme~\ref{imp} on a
    \begin{eqnarray*}
      \|v\|&=&\sup_{R\in\Omega_V}\|R(v)\|\\
      &=&\sup_{R\in J(\Omega_W)}\|R(v)\|=\sup_{S\in\Omega_W}\|S(J(v))\|\\
      &=&\|J(v)\|,
    \end{eqnarray*}
    donc $J$ est une isom\'etrie.
  \end{dem}

\end{section}

\begin{section}{Bimodules $L^\infty$ matriciellement norm\'es}
  
  \begin{deff}
    Soit $V$ un $A-B$ bimodule. On dit que $V$ est
    \emph{$L^\infty$ matriciellement norm\'e} s'il est \'equip\'e d'une
    structure d'espace d'o\-p\'e\-ra\-teurs v\'erifiant les axiomes de Ruan:
    {\def\theequation{R\addtocounter{manev}{1}\arabic{manev}}
      \begin{eqnarray}
        \label{Ruan1}
        \|\tilde{a}\tilde{v}\tilde{b}\|_n&\le&
        \|\tilde{a}\|\cdot\|\tilde{v}\|_n\cdot
        \|\tilde{b}\|\\
        \label{Ruan2}
        \left\|\left[\begin{matrix}\tilde{v}&0\\
              0&\tilde{w}
            \end{matrix}\right]\right\|_{n+m}&=&
        \max\{\|\tilde{v}\|_n,\|\tilde{w}\|_m\},
      \end{eqnarray}
      \addtocounter{equation}{-2}}%
    pour tout $n,m\ge1$, $\tilde{a}\in\mathbb{M}_n(A)$, $\tilde{v}\in
    \mathbb{M}_n(V)$, $\tilde{w}\in\mathbb{M}_m(V)$,
    $\tilde{b}\in\mathbb{M}_n(B)$.
  \end{deff}

  \begin{rem}
    \label{rem1}
    Soit $V$ un $A-B$ bimodule et pour tout $n$ entier soit $\|\cdot\|_n$ 
    une norme sur $\mathbb{M}_n(V)$. Alors $(V,(\|\cdot\|_n)_n)$ est
    un $A-B$ bimodule $L^\infty$ matriciellement norm\'e si et seulement si
    \begin{enumerate}
    \item[1.] Pour tout entiers $n<m$ l'inclusion naturelle
      \begin{equation*}
        \mathbb{M}_n(V)\subset\mathbb{M}_m(V)
      \end{equation*}
      est une isom\'etrie.
    \item[2.] Pour tout $n$ entier, $\mathbb{M}_n(V)$ est un
      $\mathbb{M}_n(A)-\mathbb{M}_n(B)$ bimodule repr\'esentable.
    \end{enumerate}
  \end{rem}

  \begin{prop}
    \label{ppp}
    Soit $R:V\to\mathcal{B}(K,H)$ un morphisme contractif de $A-B$ bimodules,
    tel que les repr\'esentations de $A$ sur $H$, respectivement de $B$ sur
    $K$, sont localement cycliques.

    Alors pour tout $\tilde{v}\in
    \mathbb{M}_n(V)$, on a l'in\'egalit\'e suivante
    \begin{equation*}
      \label{sup}
      \|R_n(\tilde{v})\|\le\sup\{\|\sum_{k=1}^n\sum_{l=1}^na_kv_{kl}b_l\|,
      \ \|\sum_{k=1}^na_ka_k^*\|\le1,\|\sum_{l=1}^nb_l^*b_l\|\le1\}.
    \end{equation*}
  \end{prop}
  
  \begin{dem}
    Pour calculer la norme de $R_n(\tilde{v})$ il suffit de
    consid\'erer des vecteurs de la forme
    \begin{equation*}
      \tilde\eta=\oplus b_k\eta,\quad\tilde\xi=\oplus a_k\xi,
    \end{equation*}
    avec $\|\tilde{\eta}\|\le1$, $\|\tilde{\xi}\|\le1$ et $\eta$, $\xi$
    appartenant \`a des sousensembles localement cycliques de $K$ et $H$ 
    respectivement.
    Soit $a=(\sum_ka_ka_k^*)^{1/2}$ et
    $b=(\sum_kb_k^*b_k)^{1/2}$.  Par le lemme~\ref{pici}, pour tout
    $0<\varepsilon<1$ il existe 
    $c_1,c_2,\ldots,c_n\in A$, $d_1,d_2,\ldots,d_n\in B$ tels que
    $a_k=c_ka^\nep$, $b_k=d_kb^\nep$ pour $k=1,2,\ldots,n$ et
    \begin{equation*}
      \sum_kc_k^*c_k=a^\dep,\quad\sum_kd_k^*d_k=b^\dep.
    \end{equation*}
    On remarque que
    \begin{equation*}
      \|b\eta\|=\|\tilde{\eta}\|,\quad\|a\xi\|=\|\tilde{\xi}\|.
    \end{equation*}
    D'autre part
    \begin{equation*}
      (R_n(\tilde{v})\tilde{\eta}|\tilde{\xi})=
      (R(\sum_{kl}c^*_kv_{kl}d_l)b^\nep\eta|
      a^\nep\xi),
    \end{equation*}
    donc
    \begin{equation*}
      |(R_n(\tilde{v})\tilde{\eta}|\tilde{\xi})|\le      
      \|\sum_{kl}c^*_kv_{kl}d_l\|\cdot\|b^\nep\eta\|\cdot\|a^\nep\xi\|.
    \end{equation*}
    Pour obtenir le r\'esultat il suffit de passer \`a la limite pour
    $\varepsilon\to0$.
\end{dem}

\begin{deff}
  Soit $V$ un $A-B$ bimodule norm\'e. On dit qu'une structure 
  $(V,(\|\cdot\|_n)_n)$ de $A-B$ bimodule $L^\infty$ matriciellement norm\'e
  est compatible si $\|\cdot\|_1=\|\cdot\|$.
\end{deff}

Il est clair que $V$ est un $A-B$ bimodule
repr\'esentable si et seulement si on peut munir $V$ d'une structure de 
bimodule
$L^\infty$ matriciellement norm\'e
compatible. Nous allons voir que parmi
toutes les structures de $A-B$ bimodules $L^\infty$ matriciellement norm\'es
compatibles avec $V$ il existe une structure minimale et une structure 
maximale (qu'on va caract\'eriser).

\begin{thm}
  \label{norm}
  Soit $V$ un $A-B$ bimodule repr\'esentable.
  Pour chaque entier $n\ge1$ on d\'efinit sur $\mathbb{M}_n(V)$ 
  \begin{equation}
    N_n(\tilde{v})=\sup\{\|\sum_{k=1}^n\sum_{l=1}^na_kv_{kl}b_l\|\ :
    \ \|\sum_{k=1}^na_ka_k^*\|\le1,\|\sum_{l=1}^nb_l^*b_l\|\le1\}.
  \end{equation}
  Alors $N_n$ est une norme pour tout $n$, $N_1(\cdot)=\|\cdot\|$ et
  $(V,(N_n)_n)$ est un $A$--$B$ bimodule $L^\infty$ matriciellement norm\'e. De
  plus, si $(V,(N'_n)_n)$ est un $A$--$B$ bimodule $L^\infty$ matriciellement
  norm\'e avec $N'_1(\cdot)= \|\cdot\|$, alors $N_n\le N'_n$, pour
  tout $n\ge1$.
\end{thm}
\begin{dem}
  On v\'erifie facilement par calcul direct que $N_n$ est une norme
  et que $N_1(v)=\|v\|$ pour tout $v\in V$.  Pour $n\ge1$ soit $\tilde{v}\in
  \mathbb{M}_n(V)$, $\tilde{\alpha}\in \mathbb{M}_n(A)$ et
  $\tilde{\beta}\in \mathbb{M}_n(B)$.  Soit $a_1,a_2,\ldots,a_n\in A$
  et $b_1,b_2,\ldots,b_n\in B$ satisfaisant \`a $\|\sum a_ka^*_k\|,\|\sum
  b^*_kb_k\|\le1$. On note $a=\left[
        \begin{matrix}
          a_1&\cdots&a_n
        \end{matrix}
      \right]\in \mathbb{M}_{1n}(A)$ et $b=\left[
          \begin{matrix}
            b_1\\
            \vdots\\
            b_n
          \end{matrix}
        \right]\in \mathbb{M}_{n1}(B)$. On aura
    \begin{eqnarray*}
      \label{submult}
      \|a.(\tilde{\alpha}\tilde{v}\tilde{\beta}).b\|&\le&
      \|a\tilde{\alpha}\|\cdot N_n(\tilde{v})\cdot\|\tilde{\beta}b\|\\
      &\le&\|\tilde{\alpha}\|\cdot N_n(\tilde{v})\cdot\|\tilde{\beta}\|,
    \end{eqnarray*}
    d'o\`u
    \begin{equation}
      \label{first_axiom}
      N_n(\tilde{\alpha}\tilde{v}\tilde{\beta})\le\|\tilde{\alpha}\|\cdot
      N_n(\tilde{v})\cdot\|\tilde{\beta}\|.
    \end{equation}
    Soient $p,q\ge1$, $\tilde{v}\in \mathbb{M}_p(V)$ et $\tilde{w}\in
    \mathbb{M}_q(V)$.  On note comme d'habitude
    $\tilde{v}\oplus\tilde{w}= \left[
      \begin{matrix}
        \tilde{v}&0\\
        0&\tilde{w}
      \end{matrix}\right]\in \mathbb{M}_{p+q}(V)$.
    D'abord on voit facilement que
    \begin{equation}
      \label{second_axiom_1}
      \max\{N_p(\tilde{v}),N_q(\tilde{w})\}\le N_{p+q}(\tilde{v}
      \oplus\tilde{w}).      
    \end{equation}
    Soit $a_1,\ldots,a_p,\ldots,a_{p+q}\in A$, $b_1,\ldots,b_p,\ldots,
    b_{p+q}\in B$ satisfaisant \`a
    \begin{equation*}
      \|\sum a_ka^*_k\|,\|\sum b^*_kb_k\|\le1.
    \end{equation*}
    Posons
    \begin{eqnarray*}
      a^1=\left[
        \begin{matrix}
          a_1&\cdots&a_p
        \end{matrix}\right]&\quad&a^2=\left[
        \begin{matrix}
          a_{p+1}&\cdots&a_{p+q}
        \end{matrix}\right],\\
      b^1=\left[
        \begin{matrix}
          b_1\\
          \vdots\\
          b_p
        \end{matrix}\right]&\quad&
      b^2=\left[
        \begin{matrix}
          b_{p+1}\\
          \vdots\\
          b_{p+q}
        \end{matrix}\right],
    \end{eqnarray*}
    Soit $u=\left[
        \begin{matrix}
          a^1&a^2
        \end{matrix}\right].(\tilde{v}\oplus\tilde{w}).\left[
        \begin{matrix}
          b^1\\
          b^2
        \end{matrix}\right]=a^1.\tilde{v}.b^1+a^2.\tilde{w}.b^2$.
      On doit d\'emontrer que
    \begin{equation}
      \label{sec}
      \|a^1.\tilde{v}.b^1+a^2.\tilde{w}.b^2\|\le\max\{N_p(\tilde{v}),
      N_q(\tilde{w})\}
    \end{equation}
    D'apr\`es le th\'eor\`eme \ref{imp}, pour tout $\varepsilon>0$ il
    existe $R:V\to B(\kk{B''},\kk{A''})$, un morphisme contractif de
    $A$--$B$ bimodules, avec $\|R(u)\|>\|u\|-\varepsilon$. Alors
      \begin{eqnarray*}
        \label{second_axiom}
        \|u\|-\varepsilon&<&\|R(a^1.\tilde{v}.b^1+a^2\tilde{w}.b^2)\|\\
        &=&\|\pi(a^1)R_p(\tilde{v})\rho(b^1)+\pi(a^2)R_q(\tilde{w})
        \rho(b^2)\|\\
        &=&\left\|\left[
            \begin{matrix}
              \pi(a^1)&\pi(a^2)
            \end{matrix}\right]\left[
            \begin{matrix}
              R_p(\tilde{v})&0\\
              0&R_q(\tilde{w})
            \end{matrix}\right]\left[
            \begin{matrix}
              \rho(b^1)\\
              \rho(b^2)
            \end{matrix}\right]\right\|\\
        &\le&\max\{\|R_p(\tilde{v})\|,\|R_q(\tilde{w})\|\}\\
        &\le&\max\{N_p(\tilde{v}),N_q(\tilde{w})\}.
      \end{eqnarray*}

      La minimalit\'e de $(V,(\|\cdot\|_n)_n)$ est imm\'ediate. 
      Soit $(V,(N'_n)_n)$ un $A-B$ bimodule $L^\infty$
      matriciellement norm\'e compatible avec $V$.
      Soient $a_1,\ldots,a_n\in A$,
      $b_1,\ldots,b_n\in B$ telles que $\|\sum_ka_ka_k^*\|\le1$,
      $\|\sum_kb_k^*b_k\|\le1$ et soit $[v_{kl}]_{kl}\in\mathbb{M}_n(V)$. 
      Alors on aura
      \begin{equation*}
        \|\sum_{kl}a_kv_{kl}b_l\|=\left\|\left[
          \begin{matrix}
            a_1\\
            \vdots\\
            a_n
          \end{matrix}\right][v_{kl}]\left[
            \begin{matrix}
              b_1\cdots b_n
            \end{matrix}\right]\right\|\le N'_n([v_{kl}]).
      \end{equation*}
    \end{dem}
    
    Les propositions~\ref{ppp} et~\ref{norm} ont comme cons\'equence 
    imm\'ediate

    \begin{cor}
      Soit $V$ un $A-B$ bimodule repr\'esentable. Alors
      \begin{eqnarray*}
        N_n(\tilde{v})&=&\sup\{\|\sum_{k=1}^n\sum_{l=1}^na_kv_{kl}b_l\|\ :
      \ \|\sum_{k=1}^na_ka_k^*\|\le1,\|\sum_{l=1}^nb_l^*b_l\|\le1\}\\
      &=&\sup_{R\in\Omega_V}\|[R(v_{kl})]\|.
    \end{eqnarray*}
  \end{cor}

    La structure minimale de $A-B$ bimodule matriciellement norm\'e
    est caract\'eris\'ee par le r\'esultat suivant:

\begin{thm}
  \label{L_infty}
  Soit $(V,(\|\cdot\|_n)_n)$ un $A-B$ bimodule $L^\infty$
  matriciellement norm\'e.
  Les affirmations suivantes sont equivalentes:
  \begin{enumerate}
  \item[(a)] La structure matricielle de $V$ est minimale.
  \item[(b)] Il existe $\Omega$ un compact et un morphisme de $A-B$
    bimodules compl\`etement isom\'etrique $V\to
    C(\Omega,\mathcal{B}(\kk{B''}, \kk{A''}))$.
  \item[(c)] Tout morphisme isom\'etrique $j:V\to
    C(\Omega,\mathcal{B}(\kk{B''},\kk{A''}))$ est compl\`etement
    isom\'etrique.
  \item[(d)] Pour $[v_{kl}]\in\mathbb{M}_n(V)$ on a
    \begin{equation*}
    \|[v_{kl}]\|_n=\sup\{\|\sum_{k=1}^n\sum_{l=1}^na_kv_{kl}b_l\|\ :
    \ \|\sum_{k=1}^na_ka_k^*\|\le1,\|\sum_{l=1}^nb_l^*b_l\|\le1\}.      
    \end{equation*}
  \item[(e)] Pour tout $A-B$ bimodule $L^\infty$ matriciellement norm\'e, 
    tout morphisme contractif de $A-B$ bimodules $W\to V$ est
    compl\`etement contractif.
  \end{enumerate}
\end{thm}
\begin{dem}
  L'\'equivalence $(a)\Leftrightarrow(d)$ a \'et\'e d\'emontr\'ee dans la
  proposition~\ref{norm}.

  $(a)\Rightarrow(c)$. Pour chaque entier $n\le1$ d\'efinissons
  \begin{equation*}
    N'_n(\tilde{v})=\|j_n(\tilde{v})\|,\quad\forall\tilde{v}\in\mathbb{M}_n(V).
  \end{equation*}
  Alors $(V,(N'_n)_n$ est un espace $L^\infty$ matriciellement norm\'e
  compatible avec $V$. Pour $\omega\in\Omega$, si on pose 
  $R_\omega(v)=j(v)(\omega)$, alors $R_\omega\in\Omega_V$. On a
  \begin{equation*}
    N'_n(\tilde{v})=\sup_{\omega\in\Omega}\|[j(v_{ij})(\omega)]\|=
    \sup_{\omega\in\Omega}\|[R_\omega(v_{ij})]\|\le N_n(\tilde{v}),
  \end{equation*}
  donc on a, pour tout entier $n$, l'\'egalit\'e $N'_n=N_n$, c'est \`a 
  dire $j$ est compl\`etement isom\'etrique.

  $(c)\Rightarrow(b)$. Cette implication est imm\'ediate, car la 
  repr\'esentation standard est une isom\'etrie, donc par hypoth\`ese elle doit \^etre compl\`etement isom\'etrique.

  $(b)\Rightarrow(a)$. Si $j$ est un morphisme compl\`etement isom\'etrique
  de $V$ dans l'espace $C(\Omega,\mathcal{B}(\kk{B''},\kk{A''}))$, alors pour
  tout $\tilde{v}\in\mathbb{M}_n(V)$
  \begin{equation*}
    \|\tilde{v}\|_n=\|j_n(\tilde{v})\|=\sup_{\omega\in\Omega}
    \|[j(v_{ij})(\omega)]\|.
  \end{equation*}
  Pour chaque $\omega\in\Omega$ soit $R_\omega(v)=j(v)(\omega)$. Alors
  $R_\omega\in\Omega_V$, donc
  \begin{equation*}
    \|\tilde{v}\|_n=\sup_{\omega\in\Omega}\|[R_\omega(v_{ij})]\|\le
    N_n(\tilde{v}).
  \end{equation*}

  $(a)\Rightarrow(e)$. Soit $W$ un $A-B$ bimodule $L^\infty$ matriciellement
  norm\'e et soit $\Phi:W\to V$ un morphisme contractif.
  Alors, pour $\tilde{w}\in\mathbb{M}_n(W)$ on aura
  \begin{eqnarray*}
    N_n(\Phi_n(\tilde{w}))&=&\sup_{R\in\Omega_V}\|[R(\Phi(w_{ij}))]\|=
    \sup_{R\in\Omega_V}\|[(R\circ\Phi)(w_{ij})]\|\\
    &\le&\sup_{S\in\Omega_W}\|[S(w_{ij})]\|,
  \end{eqnarray*}
  et on finit en utilisant la proposition~\ref{ppp} appliqu\'ee \`a $W$.

  $(e)\Rightarrow(a)$. Soit $(V,(N'_n)_n)$ un espace matriciellement norm\'e
  compatible avec $V$. Alors l'identit\'e $\mathrm{Id}_V:V\to V$ 
  est un morphisme
  contractif, et par hypoth\`ese il doit \^etre compl\`etement contractif.
  Donc $(V,(\|\cdot\|_n)_n)$ est la structure minimale de $A-B$ bimodule 
  $L^\infty$ matriciellement norm\'e.
\end{dem}

\begin{rem}
  \label{rem2}
  Soit $(V,(\|\cdot\|_n)_n)$ un $A-B$ bimodule $L^\infty$ 
  matriciellement norm\'e. Soit $k$ un entier. Alors
  $(\mathbb{M}_k(V),(\|\cdot\|_{kn})_{n\in\mathbb{N}})$ est un
  $\mathbb{M}_k(A)-\mathbb{M}_k(B)$ bimodule $L^\infty$ 
  matriciellement norm\'e.
\end{rem}

Comme cons\'equence importante des r\'esultats pr\'e\-c\'e\-dents, on peut
retrouver
le th\'eor\`eme de repr\'esentation de $A-B$ bimodules $L^\infty$ 
matriciellement norm\'es.

  \begin{thm}[\cite{ces,er1}]
    \label{effros_ruan_1}
    Soit $V$ un $A-B$ bimodule \'equip\'e d'une structure d'espace
    d'op\'erateurs. Alors $V$ v\'erifie les axiomes de Ruan $(R1)$, $(R2)$ 
    si et
    seulement si il existe un espace de Hilbert $\mathfrak{H}$, deux
    repr\'esentations $\pi:A\to\mathcal{B}(\mathfrak{H})$,
    $\rho:B\to\mathcal{B}(\mathfrak{H})$ et un morphisme de $A-B$
    bimodules compl\`etement isom\'etrique
    $J:V\to\mathcal{B}(\mathfrak{H})$ tels que
    \begin{equation*}
      J(avb)=\pi(a)J(v)\rho(b),\quad\forall a\in A,v\in V,b\in B.
    \end{equation*}
    De plus, si $A=B$ on peut choisir $\pi=\rho$.
  \end{thm}

  \begin{dem}
    Pour chaque entier $n$ soit $\mathcal{X}_n=\Omega_{\mathbb{M}_n(V)}$.
    E\-vi\-dem\-ment
    \begin{equation*}
      \mathcal{X}_n=End_{\mathbb{M}_n(A),\mathbb{M}_n(B)}
      (\mathbb{M}_n(V),\mathcal{B}(\mathbb{C}^{n^2}\otimes\kk{B''},
      \mathbb{C}^{n^2}\otimes\kk{A''}))_1,
    \end{equation*}
    car $\mathbb{C}^{n^2}$ est l'espace de la forme standard de $\mathbb{M}_n$.
    
    Consid\'erons l'application
    \begin{equation*}
      j_n:V\to\mathcal{B}(K_n,H_n)
    \end{equation*}
    donn\'ee par
    \begin{equation*}
      j_n(v)=\bigoplus_{R\in\mathcal{X}_n}R(v\oplus0),
    \end{equation*}
    o\`u
    \begin{eqnarray*}
      K_n&=&\bigoplus_{\mathcal{X}_n}(\mathbb{C}^{n^2}\otimes\kk{B''})\quad
      H_n=\bigoplus_{\mathcal{X}_n}(\mathbb{C}^{n^2}\otimes\kk{A''})\\
      v\oplus0&=&\left[
        \begin{matrix}
          v&0\\
          0&0
        \end{matrix}\right]\in\mathbb{M}_n(V),\quad\forall v\in V.
    \end{eqnarray*}
    Evidemment $j_n$ est un morphisme de $A-B$ bimodules. 
    Soit $R\in\mathcal{X}_n$. Alors pour tout $m$ entier et 
    $[v_{ij}]\in\mathbb{M}_m(V)$, par le th\'eor\`eme~\ref{L_infty} et
    la remarque~\ref{rem2} on a
    \begin{equation*}
      \|[R(v_{ij}\oplus0)]\|\le\|[v_{ij}\oplus0]\|_{mn}=\|[v_{ij}]\|_m,
    \end{equation*}
    donc $j_n$ est un morphisme compl\`etement contractif.
    D'autre part, $j_n$ est $n$-isom\'etrique, car pour tout 
    $R\in\mathcal{X}_n$ et $[v_{ij}]\in\mathbb{M}_n(V)$ on a
    \begin{equation*}
      \|[R(v_{ij}\oplus0)]\|=\|R([v_{ij}])\|.
    \end{equation*}
    En fait, $j_n$ est $k$-isom\'etrique, pour $k=1,2,\ldots,n$.

    Soit 
    \begin{equation*}
      j=\oplus_n j_n.
    \end{equation*}
    Alors $j$ est une isom\'etrie compl\`ete de $V$ \`a valeurs dans
    un espace $\mathcal{B}(K,H)$. 
    Soit $\pi_n(a)=1_{n^2}\otimes a$ et $\rho_n(b)=1_{n^2}\otimes b$.
    Evidemment 
    \begin{equation*}
      j(avb)=\pi(a)j(v)\rho(b),\quad\forall a\in A,b\in B, v\in V,
    \end{equation*}
    o\`u 
    \begin{eqnarray*}
      \pi&=&\oplus_n\oplus_{\mathcal{X}_n}\pi_n\\
      \rho&=&\oplus_n\oplus_{\mathcal{X}_n}\rho_n.
    \end{eqnarray*}
    Si $A=B$ il ne reste rien \`a d\'emontr\'er, car dans ce cas notre
    construction est telle que $\pi=\rho$. Si $A\not=B$, en consid\'erant
    l'espace $K\oplus H$, soit $J:V\to\mathcal{B}(K\oplus H)$ d\'efinie par
    \begin{equation*}
      J(v)=\left[
      \begin{matrix}
        0&0\\j(v)&0
      \end{matrix}\right]
    \end{equation*}
    et les repr\'esentations $0\oplus\pi$ respectivement $\rho\oplus0$.
    Alors $J$ est un morphisme de $A-B$ bimodules compl\`etement isom\'etrique,
    v\'erifiant les conditions de l'\'enonc\'e.    
  \end{dem}

\end{section}

\begin{section}{Bimodules normaux}

\begin{deff}
  Un \emph{$A-B$ bimodule repr\'esentable dual} est un $A-B$ bimodule
  repr\'esentable tel que $V$ est le dual d'un espace de Banach $V_*$
  et pour tout $a\in A$, $b\in B$ les applications $v\mapsto av$ et
  $v\mapsto vb$ sont $\sigma(V,V_*)$ continues.
  
  Un $A-B$ bimodule repr\'esentable dual $V$ est dit \emph{normal \`a
  gauche} (resp. \`a droite) si $A$ (resp. $B$) est une
  $W^*$-alg\`ebre et pour tout $F\in V_*$,
  l'application $a\mapsto F(av)$
  (resp.  $b\mapsto F(vb)$) appartient \`a $A_*$ (resp. $B_*$). Si $V$
  est normal \`a gauche et \`a droite on dira qu'il est \emph{normal}.
\end{deff}

Dans ce cas le lemme~\ref{pici} peut \^etre renforc\'e.

\begin{lemme}
  \label{pici_normal}
  Soit $V$ un $A-B$ bimodule normal \`a gauche (resp. \`a droite).
  Soit $F\in V^*$ avec $\|F\|\le1$ une fonctionnelle $\sigma(V,V_*)$
  continue.  Alors il existe $\phi\in S(A)$, $\psi\in S(B)$ deux
  \'etats de $A$ et $B$ respectivement tels que
    \begin{equation*}
      |F(avb)|\le\phi(aa^*)^{1/2}\psi(b^*b)^{1/2}\|v\|,
    \end{equation*}
    pour tout $a\in A$, $b\in B$ et $v\in V$, et tels que $\phi$
    (resp. $\psi$) est normal.
\end{lemme}
\begin{dem}
  
  Supposons que $V$ est normal \`a gauche et $F\in V_*$ est de norme
  $1$. Par le lemme ~\ref{pici} il existe deux \'etats $\phi$,
  $\psi$ de $A$ et $B$ respectivement tels que
    \begin{equation*}
      |F(avb)|\le\phi(aa^*)^{1/2}\psi(b^*b)^{1/2}\|v\|.
    \end{equation*}
    Par compacit\'e de $V_1$ (la boule unit\'e de $V$), il existe
    $v_0\in V_1$ de norme $1$ tel que $F(v_0)=\|F\|=1$.  
    Soit $\tau(a)=F(av_0)$.
    Alors $\tau\in A_*$ et en plus $\|\tau\|=1=\tau(1)$, donc $\tau$
    est un \'etat normal de $A$. On a aussi
    \begin{equation*}
      \tau(a)\le\phi(a)^{1/2},\quad\forall a\in A^+.
    \end{equation*}
    Pour prouver la normalit\'e de $\phi$ on utilise un crit\`ere connu (voir
    par exemple~\cite[Corollary III.3.11]{tak}).
    Consid\'erons une famille
    $(p_\alpha)_{\alpha\in I}$ de projecteurs deux \`a deux
    orthogonaux dans $A$ v\'erifiant $\sum_\alpha p_\alpha=1$. Si
    $F\subset I$ est un ensemble fini
  \begin{equation*}
    \sum_{\alpha\in F}\phi(p_\alpha)=
    \phi(\sum_{\alpha\in F}p_\alpha)\ge
    \tau(\sum_{\alpha\in F}p_\alpha)^2,
  \end{equation*}
  donc la normalit\'e de $\tau$ implique
  \begin{equation*}
    \sum_{\alpha\in I}\phi(p_\alpha)\ge1,
  \end{equation*}
  et donc $\phi$ est un \'etat normal.
\end{dem}

\begin{prop}
  \label{p2_normal}
  Soit $V$ un $A-B$ bimodule repr\'esentable dual 
  et $F\in V_*$ avec
  $\|F\|\le1$. Alors il existe $\eta\in\kk{B''}$, $\xi\in\kk{A''}$
  deux vecteurs de norme $1$ et $R:V\to B(\kk{B''},\kk{A''})$ un morphisme
  contractif de $A-B$ bimodules continu pour les topologies $*$-faibles tel que
  \begin{equation*}
    \label{main_only_dual}
    F(\cdot)=(R(\cdot)\eta|\xi).
  \end{equation*}
  Si, en outre, $V$ est normal \`a gauche (respectivement \`a droite), il
  existe $\eta\in\kk{B''}$ (resp.  $\eta\in\kk{B}$), $\xi\in\kk{A}$
  (resp. $\xi\in\kk{A''}$) de norme $1$ et $R:V\to B(\kk{B''},\kk{A})$
  (resp. $R:V\to B(\kk{B},\kk{A''})$) un morphisme contractif
  de $A$--$B$ bimodules continu pour les topologies $*$-faibles tels que
  \begin{equation*}
    \label{main1}
    F(\cdot)=(R(\cdot)\eta|\xi).
  \end{equation*}
\end{prop}
\begin{dem}
  La d\'emonstration est analogue \`a celle de la
  proposition~\ref{p2}.
  
  Supposons que $V$ est un $A-B$ bimodule repr\'esentable dual normal
  \`a gauche.  En tenant compte du lemme~\ref{pici_normal}, il
  existe un \'etat normal $\phi$ de $A$ et un \'etat $\psi$ de $B$
  tels que
  \begin{equation*}
    |F(avb)|\le\phi(aa^*)^{1/2}\psi(b^*b)^{1/2}\|v\|,
  \end{equation*}
  et en consid\'erant les repr\'esentations standard de $A$ et
  de $B''$ respectivement , $\phi$ est un \'etat vectoriel sur $\kk{A}$
  et $\psi^{**}$ est un \'etat vectoriel sur $\kk{B''}$. Soit
  $\phi=\omega_\xi$ et $\psi=\omega_\eta$. Donc il existe pour chaque
  $v$ un op\'erateur
  $R(v)\in\mathcal{B}(\overline{B''\eta},\overline{A\xi})$ tel que
  \begin{equation*}
    F(avb)=(R(v)b\eta,a^*\xi),\quad\forall a\in A,b\in B.
  \end{equation*}
  Evidemment $\overline{B''\eta}=\overline{B\eta}$.  Il est
  \'evident que $R$ est un morphisme contractif de $A-B$ bimodules.
  Il reste \`a demontrer que $R$ est continu pour les topologies $*$-faibles,
  c'est \`a dire que si
  $v_\alpha\to0$ dans la topologie $\sigma(V,V_*)$ alors
  $R(v_\alpha)\to0$ dans la topologie
  $\sigma(\mathcal{B}(\overline{B\eta},\overline{A\xi}),
  \mathcal{B}(\overline{B\eta},\overline{A\xi})_*)$. Ceci est
  \'equivalent \`a d\'emontrer que
  \begin{equation*}
    R'(\mathcal{B}(\overline{B\eta},\overline{A\xi})_*)\subset V_*,
  \end{equation*}
  o\`u $R'$ est la transpos\'ee de $R$.  Mais
  $\mathcal{B}(\overline{B''\eta},\overline{A\xi})_*=\mathcal{C}^1
  (\overline{A\xi},\overline{B\eta})$ (les op\'erateurs \`a trace de
  $\overline{A\xi}$ \`a valeurs dans $\overline{B\eta}$), 
  la dualit\'e \'etant donn\'ee par la trace.
  Comme  les op\'erateurs de rang $1$ de la forme $b\eta\otimes a\xi$
  forment un ensemble total dans
  $\mathcal{C}^1(\overline{A\xi},\overline{B\eta})$, il suffit 
  de d\'emontrer que pour tout $a\in A$ et $b\in B$, la forme lin\'eaire
  $R'(b\eta\otimes a\xi)$ est $\sigma(V,V_*)$ continue. Or, pour
  tout $v\in V$ on a
  \begin{eqnarray*}
    \left<v,R'(b\eta\otimes a\xi)\right>&=&
    \left<b\eta\otimes a\xi,R(v)\right>\\
    &=&\mathrm{Tr}(R(v)(b\eta\otimes a\xi))=
    (R(v)b\eta|a\xi)\\
    &=&F(a^*vb),
  \end{eqnarray*}
  et par hypoth\`ese les applications $F$ et $v\mapsto a^*vb$ sont
  $\sigma(V,V_*)$ continues.  Pour achever la d\'emonstration il
  suffit de remplacer $R$ par l'application $UR(\cdot)V^*$, o\`u
  $U$ et $V$ sont les deux inclusions
  $\overline{B\eta}\subset\kk{B''}$ respectivement
  $\overline{A\xi}\subset\kk{A}$.
\end{dem}

Pour simplifier, \'etant donn\'e $V$ et $W$ deux $A-B$ bimodules
duaux, on va noter avec $End_{A,B}^w(V,W)_1$ l'ensemble des morphismes
contractifs de $V$ \`a valeurs dans $W$, continus pour les topologies 
$*$-faibles.

\begin{thm}
  \label{bullshit}
  Soit $V$ un $A-B$ bimodule repr\'esentable dual. Alors
  \begin{equation*}
    \|v\|=\sup\{\|R(v)\|\quad:\quad R\in End^w_{A,B}(V,\mathcal{B}(\kk{B''},
    \kk{A''}))_1.
  \end{equation*}
  Si $V$ est normal \`a droite, alors
  \begin{equation*}
    \|v\|=\sup\{\|R(v)\|\quad:\quad R\in End^w_{A,B}(V,\mathcal{B}(\kk{B},
    \kk{A''}))_1.
  \end{equation*}
  Si $V$ est normal \`a gauche, alors
  \begin{equation*}
    \|v\|=\sup\{\|R(v)\|\quad:\quad R\in End^w_{A,B}(V,\mathcal{B}(\kk{B''},
    \kk{A}))_1.
  \end{equation*}
  Si $V$ est normal, alors
  \begin{equation*}
    \|v\|=\sup\{\|R(v)\|\quad:\quad R\in End^w_{A,B}(V,\mathcal{B}(\kk{B},
    \kk{A}))_1.
  \end{equation*}
\end{thm}
\begin{dem}
  La d\'emonstration est une application directe de la proposition
  pr\'ec\'edente.
      
  Soit $v\in V$. Supposons que $V$ est normal \`a gauche, par exemple. 
  Alors il
  existe une suite $(F_n)_n\subset V_*$ avec $\|F_n\|\le1$ telle que
  $\lim_nF_n(v)=\|v\|$. Par la proposition~\ref{p2_normal} il existe
  \begin{equation*}
  R_n\in End^w_{A,B}(V,\mathcal{B}(\kk{B''},\kk{A}))_1  
  \end{equation*}
  et des vecteurs 
  $\eta_n\in\kk{B''}$, $\xi_n\in\kk{A}$ de norme $1$ tels que
  \begin{equation*}
    F_n(v)=(R_n(v)\eta_n|\xi_n),
  \end{equation*}
  donc $\lim_n\|R_n(v)\|=\|v\|$.
\end{dem}

Avec, en plus, des hypoth\`eses de normalit\'e, le
th\'eor\`eme~\ref{repr_bimodules} devient
\begin{thm}
  \label{repr_normaux}
  Soit $V$ un $A-B$ bimodule dual. Alors il existe un espace de
  Hilbert $H$, deux repr\'esentations $\pi:A\to\mathcal{B}(H)$,
  $\rho:B\to\mathcal{B}(H)$ et un morphisme isom\'etrique de $A-B$ bimodules
  $J:V\to\mathcal{B}(H)$ continu pour les topologies
  $*$-faibles, telles que
  \begin{equation*}
    J(avb)=\pi(a)J(v)\rho(b),\quad\forall a\in A, b\in B,v\in V.
  \end{equation*}
  Si $V$ est normal \`a gauche (\`a droite), alors on peut supposer que
  $\pi$ (respectivement $\rho$) est une repr\'esentation normale
\end{thm}
\begin{dem}
  Supposons par exemple que $V$ est normal \`a gauche.  Soit
  $\mathcal{M}$ l'ensemble des morphismes (de $A-B$ bimodules) de $V$
  \`a valeurs dans $\mathcal{B}(\kk{B''},\kk{A})$, contractifs et
  continus pour les topologies $*$-faibles. Avec les notations 
  pr\'ec\'edentes, 
  \begin{equation*}
    \mathcal{M}=End^w_{A,B}(V,\mathcal{B}(\kk{B''},\kk{A}))_1.
  \end{equation*}
  Soit $\mathfrak{m}=\mathrm{card}(\mathcal{M})$ et soit
  $K=l^2(\mathfrak{m})$ On d\'efinit
  \begin{equation*}
    J:V\to\mathcal{B}(K\otimes\kk{B''},K\otimes\kk{A}),
  \end{equation*}
  par
  \begin{equation*}
    J(v)=\oplus_{R\in\mathcal{M}}R(v).
  \end{equation*}
  $J$ est une isom\'etrie par le th\'eor\`eme~\ref{bullshit}. Soient
  $\pi$ et $\rho$ les repr\'esentations de $A$ et $B$ dans $\kk{A}$
  respectivement $\kk{B''}$. Alors
  \begin{equation*}
    J(avb)=(1\otimes\pi)(a)J(v)(1\otimes\rho)(b),\quad
    \forall a\in A, b\in B, v\in V.
  \end{equation*}
  Le morphisme $J$ est \'evidemment continu pour les topologies
  $*$-faibles, \'etant une
  somme directe. Aussi, $1\otimes\pi$ est une repr\'esentation
  normale de $A$, car $\pi$ est normale. Le r\'esultat s'obtient en utilisant 
  le m\^eme argument que dans la d\'emonstration du
  th\'eor\`eme~\ref{repr_bimodules}.
\end{dem}

Enfin, comme cons\'equence directe, on donne une l\'eg\`ere
g\'en\'eralization d'un autre th\'eor\`eme de repr\'esentation
de~\cite{er1}.

\begin{deff}
  \label{dual_l_infty}
  Soit $V$ un $A-B$ bimodule $L^\infty$ matriciellement norm\'e. On dit que
  $V$ est un $A-B$ bimodule dual $L^\infty$ matriciellement norm\'e si
  il existe $(V_*,(\|\cdot\|)_n)$ un espace $L^1$ matriciellement norm\'e
  tel que $\mathbb{M}_n(V_*)^*=\mathbb{M}_n(V)$, pour tout entier $n\le1$.
\end{deff}

On rappelle qu'un espace vectoriel $V$ muni d'une syst\`eme de normes
matricielles $(\|\cdot\|_n)$ est dit $L^1$ matriciellement norm\'e s'il
v\'erifie l'axiome $(R1)$ de Ruan et si
\begin{equation*}
  \left\|\left[
    \begin{matrix}
      \tilde{v}&0\\
      0&\tilde{w}
    \end{matrix}\right]\right\|=\|\tilde{v}\|_n+\|\tilde{w}\|_m,
\end{equation*}
pour tous les entiers $n$, $m$ et pour tout $\tilde{v}\in\mathbb{M}_n(V)$,
$\tilde{w}\in\mathbb{M}_m(V)$.

Pour un espace d'op\'erateurs donn\'e $V$, il est possible que $V$ poss\`ede
un pr\'edual $V_*$ tel que la boule unit\'e de $\mathbb{M}_2(V)$ n'est
pas ferm\'ee pour la topologie g\'en\'er\'ee par les seminormes de la forme
\begin{equation*}
  \left[\begin{matrix}
    v_{11}&v_{12}\\
    v_{21}&v_{22}
  \end{matrix}\right]\mapsto
\left\|\left[\begin{matrix}
    f(v_{11})&f(v_{12})\\
    f(v_{21})&f(v_{22})
\end{matrix}\right]\right\|,
\end{equation*}
o\`u $f\in V_*$.

\begin{rem}
  Si $V$ est un $A-B$ bimodule dual $L^\infty$ matriciellement norm\'e
  alors pour tout entier $n$, le $\mathbb{M}_n(A)-
  \mathbb{M}_n(B)$ bimodule $\mathbb{M}_n(V)$ est dual. Si $V$ est normal
  \`a gauche, alors $\mathbb{M}_n(V)$ est normal \`a gauche.
  La m\^eme remarque est valable pour le cas normal \`a droite.
\end{rem}
\begin{thm}
  \label{effros_ruan_2}
  Soit $V$ un $A-B$ bimodule dual $L^\infty$ matriciellement norm\'e.
  Alors il existe
  un espace de Hilbert $\mathfrak{H}$, deux repr\'esentations
  $\pi:A\to\mathcal{B}(\mathfrak{H})$,
  $\rho:B\to\mathcal{B}(\mathfrak{H})$ et un morphisme de $A-B$
  bimodules compl\`etement isom\'etrique et continu pour les topologies
  $*$-faibles
  $J:V\to\mathcal{B}(\mathfrak{H})$ tels que
  \begin{equation*}
    J(avb)=\pi(a)J(v)\rho(b),\quad\forall a\in A,v\in V,b\in B.
  \end{equation*}
  De plus, si $A=B$ on peut choisir $\pi=\rho$.  Si $V$ est normal \`a
  gauche (\`a droite), on peut choisir $\pi$ (resp. $\rho$) une
  repr\'esentation normale.
\end{thm}
\begin{dem}
  La d\'emonstration est analogue \`a celle du
  th\'e\-o\-r\`eme~\ref{effros_ruan_1}.  Par exemple, si $V$ est normal \`a
  gauche, pour chaque entier $n$ soit
  \begin{eqnarray*}
    \mathcal{X}_n=End^w_{\mathbb{M}_n(A),\mathbb{M}_n(B)}
    (\mathbb{M}_n(V),\mathcal{B}(\mathbb{C}^{n^2}\otimes\kk{A},
    \mathbb{C}^{n^2}\otimes\kk{B''})).
  \end{eqnarray*}
  Evidemment
  $\mathbb{C}^{n^2}$ est l'espace de la forme standard de
  $\mathbb{M}_n$.
  Consid\'erons l'application
  \begin{equation*}
    j_n:V\to\mathcal{B}(K_n,H_n)
  \end{equation*}
  donn\'ee par
  \begin{equation*}
    j_n(v)=\bigoplus_{R\in\mathcal{X}_n}R(v\oplus0),
  \end{equation*}
  o\`u
  \begin{eqnarray*}
    K_n&=&\bigoplus_{\mathcal{X}_n}(\mathbb{C}^{n^2}\otimes\kk{B''})\\
    H_n&=&\bigoplus_{\mathcal{X}_n}(\mathbb{C}^{n^2}\otimes\kk{A}).
  \end{eqnarray*}
  Soit $\pi_n(a)=1_{n^2}\otimes a$ et $\rho_n(b)=1_{n^2}\otimes b$.
  Alors $j_n$ est un morphisme de $A-B$ bimodules et
  \begin{equation*}
    j_n(avb)=\pi'_n(a)j_n(v)\rho'_n(b),\quad\forall a\in A,b\in B,v\in V,
  \end{equation*}
  o\`u $\pi'_n$ et $\rho'_n$ sont des sommes directes des repr\'esentations
  $\pi_n$ respectivement $\rho_n$.
  En particulier, $\pi'_n$ est une repr\'esentation normale.
  Soit $R\in\mathcal{X}_n$. Alors pour tout $m$ entier et
  $[v_{ij}]\in\mathbb{M}_m(V)$, toujours en utilisant le 
  th\'eor\`eme~\ref{L_infty} et
  la remarque~\ref{rem2} on a
  \begin{equation*}
    \|[R(v_{ij}\oplus0)]\|\le\|[v_{ij}\oplus0]\|_{mn}=\|[v_{ij}]\|_m,
  \end{equation*}
  donc $j_n$ est un morphisme compl\`etement contractif.  D'autre
  part, $j_n$ est $n$-isom\'etrique, par
  le th\'eor\`eme~\ref{bullshit}, car il est facile de voir que
  $\mathbb{M}_n(V)$ est un $\mathbb{M}_n(A)-\mathbb{M}_n(B)$ bimodule
  dual. Aussi, $j_n$ est continu pour les topologies $*$-faibles, car c'est
  une somme directe de telles applications.
  
  Soit
  \begin{equation*}
    j=\oplus_n j_n.
  \end{equation*}
  Alors $j$ est une isom\'etrie compl\`ete de $V$ \`a valeurs dans
  un espace $\mathcal{B}(K,H)$. Evidemment
  \begin{equation*}
    j(avb)=\pi(a)j(v)\rho(b),\quad\forall a\in A,b\in B, v\in V,
  \end{equation*}
  o\`u
  \begin{equation*}
    \pi=\oplus_n\pi'_n,\quad
    \rho=\oplus_n\rho'_n.
  \end{equation*}
  Evidemment $\pi$ est une repr\'esentation normale de $A$.
  Si $A=B$ il ne reste rien \`a d\'emontrer, car dans ce cas notre
  construction est telle que $\pi=\rho$. Si $A\not=B$, en
  consid\'erant l'espace $K\oplus H$, soit
  $J:V\to\mathcal{B}(K\oplus H)$ d\'efinie par
  \begin{equation*}
    J(v)=\left[
      \begin{matrix}
        0&0\\j(v)&0
      \end{matrix}\right]
  \end{equation*}
  et les repr\'esentations $0\oplus\pi$ respectivement
  $\rho\oplus0$.  Alors $J$ est un morphisme de $A-B$ bimodules
  compl\`etement isom\'etrique, v\'erifiant les conditions de
  l'\'enonc\'e.
\end{dem}

\end{section}

\bibliographystyle{plain} 
\nocite{wit1}
\bibliography{standard}

\vspace{2cm}
\font\cucu=cmcsc10
\newbox\luluc
\newbox\lulu
\setbox\luluc=\vtop{\hsize=7cm\parindent=0pt\cucu Institutul de Matematic\u a
  al Academiei Rom\^ane\\
  C.P. 1--764\\
  Bucure\c sti\\
  Rom\^ania}
\setbox\lulu=\vtop{\hsize=7cm\parindent=0pt\cucu Universit\'e d'Orl\'eans\\
  U.F.R. Sciences\\
  Dept. de Math\'ematiques\\
  B.P. 6759 ORLEANS CEDEX 2\\
  FRANCE}
\begin{tabular}[hbp!]{l@{\ :\ \ }l}
Adresse&\box\luluc\\
\noalign{\vspace{0.5cm}}
E-mail&\texttt{cipop@stoilow.imar.ro}\\
\noalign{\vspace{1cm}}
Adresse actuelle&\box\lulu\\
\noalign{\vspace{0.5cm}}
E-mail&\texttt{pop@labomath.univ-orleans.fr}
\end{tabular}

\end{document}